\documentclass{article}

\usepackage{amsmath,amssymb,latexsym}
\usepackage{url}

\newtheorem{theorem}{Theorem}
\newcommand{\bt}{\begin{theorem}}
\newcommand{\et}{\end{theorem}}
\newtheorem{lemma}{Lemma}
\newcommand{\bl}{\begin{lemma}}
\newcommand{\el}{\end{lemma}}
\newtheorem{corollary}{Corollary}
\newcommand{\bc}{\begin{corollary}}
\newcommand{\ec}{\end{corollary}}
\newcommand{\beq}{\begin{equation}}
\newcommand{\eeq}{\end{equation}}
\newcommand{\benum}{\begin{enumerate}}
\newcommand{\eenum}{\end{enumerate}}
\newcommand{\N}{\ensuremath{ \mathbf N }}
\newcommand{\Z}{\ensuremath{\mathbf Z}}
\newcommand{\Q}{\ensuremath{\mathbf Q}}
\newcommand{\mca}{\ensuremath{ \mathcal A}}

\DeclareMathOperator{\card}{\text{card}}

\begin{document}

\title{On the fractional parts of roots 
of positive real numbers}
\author{Melvyn B. Nathanson}
\date{}
\maketitle

\begin{abstract}
Let  $[\theta]$ denote the integer part and $\{\theta\}$  the fractional part of the real number $\theta$.  
For $\theta > 1$ and $\{ \theta^{1/n} \} \neq 0$,  define 
 $M_{\theta} (n) = [1/\{ \theta^{1/n} \}]$. 
The arithmetic  function $M_{\theta} (n)$ is eventually increasing, 
and $\lim_{n\rightarrow \infty} M_{\theta}(n)/n = 1/\log \theta$.  
Moreover, $M_{\theta} (n)$ is ``linearly periodic'' if and only if 
$\log \theta$ is rational.  
Other results and problems concerning the function $M_{\theta}(n)$ are discussed.
\end{abstract}

\section{The sequence of roots and the arithmetic function 
$M_{\theta}(n)$}

Let \N, $\N_0$, and \Z\ denote the positive integers, nonnegative integers, and integers, respectively.  
An \emph{arithmetic function} is a function whose domain is the set \N\ of positive integers.  
Let ${\theta}$ be a real number, 
and let $[{\theta}]$ denote the integer part of  ${\theta}$ 
and  $\{{\theta}\}$  the fractional part of ${\theta}$.  
Thus, ${\theta} = [{\theta}] + \{ {\theta}\}$, where $[{\theta}] \in \Z$ and $0 \leq \{{\theta}\} < 1$.
Let $\| {\theta} \| = \min( \{{\theta}\},1-\{{\theta}\})$ denote the distance from ${\theta}$ to the nearest integer.

A famous theorem of Koksma~\cite{koks35}  
(see Kuipers and 
Niederreiter~\cite[Corollary~4.2]{kuip-nied74}) 
states that the sequence of the fractional parts of the $n$th powers of ${\theta}$, that is, 
$( \{{\theta}^n\} )_{n=1}^{\infty}$, is uniformly distributed in the interval $[0,1]$  for almost all real numbers ${\theta} > 1$.  
Nonetheless, there is no known number 
${\theta}$ whose powers are uniformly distributed modulo 1.
It is a famous unsolved problem to understand the distribution of  the  fractional parts of the powers of a rational number, and, in particular, of $3/2$ 
(cf.~\cite{bake-coat75,dubi09,flat-laga-poll95, habs03, nath72a, wags79, zudi07}).

There is a large body of research on the fractional parts of powers, 
but there seems to have been no investigation of the dual problem of the distribution of the fractional parts of the $n$th roots of a positive real number ${\theta} \neq 1$.   
Mahler and Szekeres~\cite{mahl-szek67} and Bugeaud and Dubickas~\cite{buge-dubi08} have considered the distribution modulo 1 of the sequence $\left( \| {\theta}^n\|^{1/n} \right)_{n=1}^{\infty}$, but this is different from the sequences that will be considered in this paper.

Let $\theta$ be a  positive real number.
For every positive integer $n$ such that $\{ \theta^{1/n} \} \neq 0$, we define the arithmetic function  
\[
M_{\theta}(n) = 
 \left[ \frac{1}{  \{ {\theta}^{1/n} \} } \right]. 
\] 
Let $M_{\theta}(n) =  \infty$ if $\{ \theta^{1/n} \} = 0$.  
Note that  $M_{\theta}(n) = \infty$ for infinitely many $n \in \N$ if and only if $\theta = 1$. 
We observe  that if $0 <   \theta < 1$ and $n >  -\log  \theta /\log 2$, 
then $1/2 <  \theta^{1/n} =  \{   \theta^{1/n} \} < 1$.  
It follows that $1 <  \left\{ \theta^{1/n}\right\}^{-1} < 2$ and 
$M_{\theta}(n) = \left[ \left\{\theta^{1/n} \right\}^{-1} \right] =1$.
Thus, the function $M_{\theta}(n)$ is eventually constant for $0 < \theta < 1$, and so it suffices to consider only $\theta > 1$.  

For $\theta > 1$, let $N_0(\theta)$ denote the smallest integer $n$ such that $n > \log\theta/\log 2$.  
If $n \geq N_0(\theta)$, 
then $1 < \theta^{1/n} < 2$ and so
\[
0 <   \{ \theta^{1/n} \} =  \theta^{1/n} - 1 < 1
\]
and
\[
M_{\theta}(n) = 
 \left[ \frac{1}{   {\theta}^{1/n} - 1 } \right]. 
\]

We can use Maple to compute the function $M_{\theta}(n)$ for various  $\theta$ and for $n$ from 1 to 90.
Here is the data for $\theta = 3/2$, $2$, $17$, and $\pi$ with
$1 \leq n \leq 90$.
We put a box around $M_{\theta}(N_0(\theta))$. 
We obtain the following eventually increasing sequences of integers.
\[
\boxed{\theta = \frac{3}{2}}
\]
\[
\begin{array}{ccccccccccccccc}  
\boxed{2} & 4 & 6 & 9 & 11 & 14 & 16 & 19 & 21 & 24 & 26 & 29 & 31 & 34 & 36 \\
38 & 41 & 43 & 46 & 48 & 51 & 53 & 56 & 58 & 61 & 63 & 66 & 68 & 71 & 73\\
75 & 78 & 80 & 83 & 85 & 88 & 90 & 93 & 95 & 98 & 100 & 103 & 105 & 108 & 110 \\
112 & 115 & 117 & 120 & 122 & 125 & 127 & 130 & 132 & 135 & 137 & 140 & 142 & 145 & 147 \\
149 & 152 & 154 & 157 & 159 & 162 & 164 & 167 & 169 & 172 & 174 & 177 & 179 & 182 & 184 \\
186 & 189 & 191 & 194 & 196 & 199 & 201 & 204 & 206 & 209 & 211 & 214 & 216 & 219 & 221
 \end{array}
\]
\\
\[
\boxed{\theta = 2}
\]
\[
\begin{array}{ccccccccccccccc}  
\infty & \boxed{2} & 3 & 5 & 6 & 8 & 9 & 11 & 12 & 13 & 15 & 16 & 18 & 19 & 21 \\
22 & 24 & 25 & 26 & 28 & 29 & 31 & 32 & 34 & 35 & 37 & 38 & 39 & 41 & 42 \\
44 & 45 & 47 & 48 & 49 & 51 & 52 & 54 & 55 & 57 & 58 & 60 & 61 & 62 & 64 \\
65 & 67 & 68 & 70 & 71 & 73 & 74 & 75 & 77 & 78 & 80 & 81 & 83 & 84 & 86 \\
87 & 88 & 90 & 91 & 93 & 94 & 96 & 97 & 99 & 100 & 101 & 103 & 104 & 106 & 107 \\
109 & 110 & 112 & 113 & 114 & 116 & 117 & 119 & 120 & 122 & 123 & 125 & 126 & 127 & 129
\end{array}
\]
\\
\[
\boxed{\theta = 17}
\]
\[
\begin{array}{ccccccccccccccc}  
\infty & 8 & 1 & 32 & \boxed{1} & 1 & 2 & 2 & 2 & 3 & 3 & 3 & 4 & 4 & 4 \\
     5 & 5 & 5 & 6 & 6 & 6 & 7 & 7 & 7 & 8 & 8 & 9 & 9 & 9 & 10 \\
10 & 10 & 11 & 11 & 11 & 12 & 12 & 12 & 13 & 13 & 13 & 14 & 14 & 15 & 15 \\
15 & 16 & 16 & 16 & 17 & 17 & 17 & 18 & 18 & 18 & 19 & 19 & 19 & 20 & 20 \\
21 & 21 & 21 & 22 & 22 & 22 & 23 & 23 & 23 & 24 & 24 & 24 & 25 & 25 & 25 \\
26 & 26 & 27 & 27 & 27 & 28 & 28 & 28 & 29 & 29 & 29 & 30 & 30 & 30 & 31
\end{array}
\]
\\
\[
\boxed{\theta = \pi}
\]
\[
\begin{array}{ccccccccccccccc}  
7 & \boxed{1} & 2 & 3 & 3 & 4 & 5 & 6 & 7 & 8 & 9 & 9 & 10 & 11 & 12 \\
13 & 14 & 15 & 16 & 16 & 17 & 18 & 19 & 20 & 21 & 22 & 23 & 23 & 24 & 25 \\
26 & 27 & 28 & 29 & 30 & 30 & 31 & 32 & 33 & 34 & 35 & 36 & 37 & 37 & 38 \\
39 & 40 & 41 & 42 & 43 & 44 & 44 & 45 & 46 & 47 & 48 & 49 & 50 & 51 & 51 \\
52 & 53 & 54 & 55 & 56 & 57 & 58 & 58 & 59 & 60 & 61 & 62 & 63 & 64 & 65 \\
65 & 66 & 67 & 68 & 69 & 70 & 71 & 72 & 72 & 73 & 74 & 75 & 76 & 77 & 78
\end{array}
\]

It is an open problem to understand and ``predict'' these sequences of integers.  The goal of this paper is to obtain  basic results about the function $M_{\theta}(n)$ and to ask some questions suggested by the experimental data.  
We shall prove that $M_{\theta}(n)$ is, for fixed $\theta$,  
an eventually increasing function of $n$ and, for fixed $n$, an eventually decreasing function of $\theta$ 
(Theorem~\ref{FPR:theorem:monotonicity}), 
that $M_{\theta}(n) \sim n/\log \theta$ 
(Theorem~\ref{FPR:theorem:AsymptoticGrowth} ), 
and that $M_{\theta}(n) = M_{\psi}(n)$ for infinitely many $n$ 
if and only if $\theta = \psi$  
(Theorem~\ref{FPR:theorem:uniqueness}).
We shall also prove that  $M_{\theta}(n)$ is strictly increasing if and only if $\theta \leq e$ (Theorem~\ref{FPR:theorem:StrictIncrease}).
For almost all $\theta > 1$, the computational data for the function $M_{\theta}(n)$ have no obvious pattern, but for certain $\theta$
we find that $M_{\theta}(n)$ is \emph{eventually linearly periodic} 
in the sense that there exist positive integers $k$, $\ell$, and $n_0$ 
such that 
\beq   \label{FPR:IntroELP}
M_{\theta}(n+k) = M_{\theta}(n) + \ell
\eeq
for all $n \geq n_0$.    
A fundamental result of this paper 
(Theorem~\ref{FPR:theorem:LinearPeriodicity}) is that identity~\eqref{FPR:IntroELP} holds if and only if 
$\theta = e^{k/\ell}$, and that there is a simple algorithm 
(Theorem~\ref{FPR:theorem:algorithm}) to compute the periodic pattern.

\section{Growth, asymptotics, and uniqueness of the function  
$M_{\theta}(n)$}    \label{FPR:section:growth}  

\bt        \label{FPR:theorem:monotonicity}
Let $\theta$ and $\psi$ be real numbers such that $1 < \psi < \theta$.  
If  $n > \log \theta/ \log 2$, then 
\[
M_{\theta}(n) \leq M_{\psi}(n)
\]
and
\[
M_{\theta}(n) \leq M_{\theta}(n+1).
\]
Moreover,
\[
M_{\theta}(n) >  \frac{n}{{\theta}-1} - 1
\]
and so
$\lim_{n\rightarrow \infty} M_{\theta}(n) = \infty.$
\et

\textit{Proof.}
The inequality $n > \log \theta/ \log 2 > \log \psi / \log 2$ 
 implies  that $\psi^{1/n}  < \theta^{1/n}$ and 
\[
M_{\theta}(n) \leq  \frac{1}{ {\theta}^{1/n}-1 } 
 < \frac{1}{ {\psi}^{1/n}-1 } 
 < M_{\psi}(n) + 1
\]
and so $M_{\theta}(n) \leq  M_{\psi}(n)$.  
Similarly, 
\[
1 < {\theta}^{1/(n+1)} < {\theta}^{1/n} \leq {\theta}  
\]
and so 
\[
M_{\theta}(n) \leq \frac{1}{ {\theta}^{1/n} -1 } < \frac{1}{ {\theta}^{1/(n+1)} -1} < M_{\theta}(n+1) + 1
\]
and  $M_{\theta}(n) \leq M_{\theta}(n+1).$

Let  $\varepsilon_n = {\theta}^{1/n}  - 1 > 0$.
Applying the binomial theorem, we obtain 
${\theta} = (1+\varepsilon_n)^n > 1 + n\varepsilon_n$.  
Equivalently,  
\[
0 < \varepsilon_n  < \frac{{\theta}-1}{n}.
\]
It follows that  
\[
M_{\theta}(n) >  \frac{1}{ {\theta}^{1/n} -1 } - 1
=  \frac{1}{\varepsilon_n } - 1
>  \frac{n}{{\theta}-1} - 1
\]
and so 
$\lim_{n\rightarrow \infty} M_{\theta}(n) = \infty.$
This completes the proof.
$\blacksquare$ 

\bt    \label{FPR:theorem:AsymptoticGrowth}
If  $\theta > 1$, then
\[
\lim_{n\rightarrow \infty} \frac{M_{\theta}(n)}{n} = \frac{1}{\log \theta}.
\]
\et

\textit{Proof.}
For all real numbers $x$ we have 
\[
\lim_{n\rightarrow \infty}  \left(1+\frac{ x}{n} \right)^n = e^x.
\]
Let $0 < \varepsilon < 1$.  For $x > 0$, the inequality 
\[
\lim_{n\rightarrow \infty}  \left(1+\frac{(1-\varepsilon) x}{n} \right)^n 
= e^{ (1- \varepsilon) x } < e^x  <  e^{ (1+ \varepsilon) x }  
= \lim_{n\rightarrow \infty} \left(1+\frac{(1+\varepsilon) x}{n} \right)^n 
\]
implies that there exists an integer 
$N(\varepsilon) > \log \theta/ \log  2$ such that 
\[
\left(1+\frac{ (1 - \varepsilon)x}{n} \right)^n < e^x 
< \left(1+\frac{(1+\varepsilon) x}{n} \right)^n 
\]
for all $n \geq N(\varepsilon)$.  
Taking $n$th roots, subtracting 1,  and reciprocating, we obtain 
\[
\frac{n}{(1+\varepsilon)x} < \frac{1}{e^{x/n} - 1} 
< \frac{n}{  (1 - \varepsilon) x}.
\]
If $\theta > 1$, then $x = \log \theta > 0$, and so 
\[
\frac{n}{(1+\varepsilon) \log \theta} < \frac{1}{\theta^{1/n} - 1} < \frac{n}{  (1 - \varepsilon)  \log \theta}.
\]
Then 
\[
\frac{n}{(1+\varepsilon) \log \theta} - 1<  M_{\theta}(n) 
< \frac{n}{  (1 - \varepsilon) \log \theta}.
\]
Equivalently, 
\[
\frac{1}{(1+\varepsilon) \log \theta} - \frac{1}{n} 
<  \frac{M_{\theta}(n)}{n} < \frac{1}{  (1 - \varepsilon) \log \theta} 
\]
for all integers $n \geq N(\varepsilon)$.  It follows  that 
\[
\frac{1}{(1+\varepsilon) \log \theta} 
\leq \liminf_{n\rightarrow \infty }  \frac{M_{\theta}(n)}{n} 
\leq \limsup_{n\rightarrow \infty } \frac{M_{\theta}(n)}{n} 
\leq  \frac{1}{  (1 - \varepsilon) \log \theta} 
\]
for all $\varepsilon > 0$, and so 
$\lim_{n\rightarrow \infty } M_{\theta}(n)/n= 1/\log \theta$.  
This completes the proof.
$\blacksquare$  

\bt       \label{FPR:theorem:uniqueness}
Let $1 < \psi  \leq \theta$ be real numbers.  
Let $(n_i)_{i=1}^{\infty}$ be a strictly increasing sequence of 
positive integers and let  $( \varepsilon_i)_{i=1}^{\infty}$ be a sequence of 
integers such that $n_i +  \varepsilon_i \geq 1$ for all $i$ and $\lim_{i\rightarrow \infty}  \varepsilon_i/n_i = 0$.  
If  $M_{\psi}(n_i +  \varepsilon_i ) = M_{\theta}(n_i)$ for all  $i$, then $\psi  = \theta $.
In particular, if  $M_{\psi}(n) = M_{\theta}(n)$ for infinitely many positive integers $n$, then $\psi  = \theta $.
\et

\textit{Proof.}
If  $M_{\psi}(n_i +  \varepsilon_i ) = M_{\theta}(n_i)$ 
for $i=1,2,\ldots,$ then
\begin{align*}
\frac{1}{\log \psi} 
& = \lim_{i \rightarrow \infty} \frac{M_{\psi}(n_i  +  \varepsilon_i)    }{n_i +  \varepsilon_i } 
= \lim_{i \rightarrow \infty} \frac{M_{\theta}(n_i)}{n_i +  \varepsilon_i}  \\
& = \lim_{i \rightarrow \infty} \frac{M_{\theta}(n_i)}{n_i}
 \lim_{i \rightarrow \infty} \frac{ n_i}{n_i +  \varepsilon_i} \\
& = \frac{1}{\log \theta}
\end{align*}
and so $\psi  = \theta $.  
In particular, if $\varepsilon_i = 0$ for all $i$, 
then  $M_{\psi}(n) = M_{\theta}(n)$ for infinitely many positive integers $n$ only if $\psi  = \theta $.
This completes the proof.  
$\blacksquare$ 

The asymptotic estimate for $M_{\theta}(n)$ given  in 
Theorem~\ref{FPR:theorem:AsymptoticGrowth}
can be sharpened for $1 < \theta \leq e$.

\bt    \label{FPR:theorem:main}
Let $\theta$ be a real number such that $1 < \theta \leq e$.  For every integer $n > \log \theta/\log  2$, 
\beq   \label{FPR:ineq3}
\frac{n-1}{\log \theta} < \frac{1}{\theta^{1/n}-1} < \frac{n}{\log \theta}   
\eeq
and
\beq   \label{FPR:ineq4}
\left[ \frac{n-1}{\log \theta} \right] \leq  M_{\theta}(n) 
< \frac{n}{\log \theta}.   
\eeq
\et

\textit{Proof.}
If  $1 < {\theta} \leq e$ and $x = \log \theta$, 
then   $0 < x \leq 1$.  
By Lemma~\ref{FPR:lemma:exp(x)}  in Appendix~\ref{FPR:appendix:exp}, for every integer $n \geq 2$ 
we have
\[
\left(1 + \frac{x}{n}\right)^n < e^x < \left(1 + \frac{x}{n-1}\right)^n 
\]
and so
\[
\frac{n-1}{x} < \frac{1}{e^{x/n}-1} < \frac{n}{x} .
\]
Equivalently, 
\[
\frac{n-1}{\log {\theta}} < \frac{1}{{\theta}^{1/n}-1} < \frac{n}{\log {\theta}} .
\]
This proves~\eqref{FPR:ineq3}, 
and inequality~\eqref{FPR:ineq3} implies~\eqref{FPR:ineq4}.
$\blacksquare$ 

\bc   \label{FPR:corollary:e}
$M_e(1) = 1 $ and
\[
M_e(n) = n-1
\] 
for every integer $n \geq 2$.
\ec

\textit{Proof.}
We have $1/2 < \{ e\} = e-2  < 1$ and so $1 < (e-2)^{-1} < 2$ and $M_e(1) = 1$.  
For $n \geq 2$ we have $1 < e^{1/n} \leq e^{1/2} <  2$.  
Applying Theorem~\ref{FPR:theorem:main} with $\theta = e$, 
$\log \theta = 1$, and $n \geq 2$, we obtain
\[
n-1\leq M_e(n) < n
\]
and so $M_e(n) = n-1$.
This completes the proof.
$\blacksquare$ 

\bt  \label{FPR:theorem:BoundedGaps}
For $x \geq 1$, define the function 
\beq   \label{FPR:g(x)}
g(x) = \frac{1}{\log\left(1 + \frac{1}{x}\right)}.
\eeq
Let $\theta > 1$ and $n > \log \theta / \log  2$.  
Then  $M_{\theta}(n) = x$ if and only if 
\beq  \label{FPR:inverse}
g(x) \log \theta  \leq n < g(x+1) \log \theta.
\eeq
If $M_{\theta}(n) = x$, then 
\beq   \label{FPR:MainIneq}
\left(x+\frac{1}{2} - \frac{1}{x}\right) \log \theta  \leq n <  \left(x+\frac{3}{2}\right) \log \theta.
\eeq
\et

\textit{Proof.}
Let $n > \log \theta / \log  2$.  
We have  $M_{\theta}(n) = x$ if and only if 
\[
x \leq  \frac{1}{ \theta^{1/n} - 1}  < x+1.
\]
Solving this equation for $n$, we obtain~\eqref{FPR:inverse}.

By Lemma~\ref{FPR:lemma:g(x)ineq} 
in Appendix~\ref{FPR:appendix:exp}, for $x \geq 1$ 
the function $g(x)$ is positive and strictly increasing, 
and satisfies  inequality~\eqref{FPR:g1}.
Inserting the estimates from~\eqref{FPR:g1} 
into~\eqref{FPR:inverse} gives~\eqref{FPR:MainIneq}.
$\blacksquare$ 

\bc    \label{FPR:corollary:BoundedGapsIneq}
Let $\theta > 1$.  
If $n_1$ and $n_2$ are integers such that 
$n_2 > n_1 > \log \theta / \log  2$, then 
\beq   \label{FPR:BoundedGapsIneq}
M_{\theta}(n_2) - M_{\theta}(n_1) < \frac{n_2 - n_1}{\log \theta} +  \frac{3}{2}.
\eeq	
\ec

\textit{Proof.}
Let $n_2 > n_1 > \log \theta / \log  2$. 
Theorem~\ref{FPR:theorem:monotonicity} implies that 
$x_1 = M_{\theta}(n_1) \leq M_{\theta}(n_2) = x_2$.  
If $x_1 < x_2$, then, by~\eqref{FPR:inverse},
\[ 
n_1 < g(x_1+1) \log \theta \leq g(x_2) \log \theta  \leq n_2.
\]
Applying~ inequality~\eqref{FPR:g2} from Appendix A, 
we obtain 
\[
\frac{n_2 - n_1}{\log \theta} 
>  g(x_2) -g(x_1+1)  
 > x_2 - x_1 - 1 - \frac{1}{x_2} 
 \geq x_2 - x_1 - \frac{3}{2}.
\]
This completes the proof.
$\blacksquare$ 

\bc     \label{FPR:corollary:Beatty-like}
Let $\theta > 1$.
If   $n > \log\theta/\log 2$, then
\[
M_{\theta}(n) = \left[ \frac{n}{\log\theta}  \pm \frac{1}{2} \right].
\]
\ec

\textit{Proof.}
Let 
\[
L_n = \left[ \frac{n}{\log\theta}  - \frac{1}{2} \right]  
\] 
and 
\[
\lambda_n = \left\{ \frac{n}{\log\theta}  - \frac{1}{2} \right\}  
=  \frac{n}{\log\theta}  - \frac{1}{2} - L_n  \in [0,1).
\]
Let $x = M_{\theta}(n)$.  
Rearranging inequality~\eqref{FPR:MainIneq}, we obtain
\[
L_n-1 \leq L_n -1 + \lambda_n <  x 
\leq  L_n + \lambda_n+ \frac{1}{x} < L_n+2.
\]
Because $x$ is an integer, we have $x = L_n$ or $x = L_n+1$.  
This completes the proof.  
$\blacksquare$ 

An arithmetic function $f(n)$ is \emph{eventually strictly increasing} if there exists an integer $n_0$ such that $f(n) < f(n+1)$ 
for all $n \geq n_0$.

\bt    \label{FPR:theorem:StrictIncrease}
Let $\theta > 1$.  
The arithmetic function $M_{\theta}(n)$ is eventually strictly increasing  if and only if $\theta \leq e$.
\et

\textit{Proof.}
By Corollary~\ref{FPR:corollary:e},  
we have $M_e(n) = n-1$ for $n \geq 2$, 
and so $M_e(n)$ is eventually strictly increasing.

Let $1 < \theta < e$.  Then $0 < \log \theta < 1$.  
By Theorem~\ref{FPR:theorem:BoundedGaps},  
for every integer $n > \log \theta / \log  2$, 
we have  $M_{\theta}(n) = x$ if and only if 
\[
g(x) \log \theta  \leq n < g(x+1) \log \theta.
\]
The length of this interval is $(g(x+1)  - g(x) ) \log \theta$.  
Applying~\eqref{FPR:g2} with  $y = x+1$, we obtain
\[
1 - \frac{1}{x+1} <  g(x+1) - g(x) < 1 + \frac{1}{x}.  
\]
Because $\lim_{n\rightarrow\infty} M_{\theta}(n) = \infty$, 
we have
\[
x = M_{\theta}(n) > \frac{ \log \theta}{1 - \log \theta} > 0 
\]
for all sufficiently large $n$, and so   
\[
0 < (g(x+1) - g(x) )\log \theta < \left( 1 + \frac{1}{x}\right) \log \theta  < 1.
\]
This implies that the interval $[g(x) \log \theta ,  g(x+1) \log \theta)$ contains at most one integer; that is, there is at most one integer $n$ such that $M_{\theta}(n) = x$.  
This means that the function 
$M_{\theta}(n)$ is eventually strictly increasing.

Let $\theta > e$.  Then  $\log  \theta > 1$.  
If $x$ and $y$ are positive integers such that 
\[
y \geq x + \frac{\log\theta + 1}{\log\theta -1}
\]
then the lower bound in~\eqref{FPR:g2}  gives 
\begin{align*}
\left( g(y) - g(x)  \right) \log \theta 
& > \left( y - x -\frac{1}{y} \right) \log \theta \\
& > (y - x - 1) \log \theta \\
& \geq y-x+1.
\end{align*}
If $n \in\N$ and $g(x)\log \theta \leq n < g(y)\log \theta$, then $x \leq M_{\theta}(n) \leq y-1$.  
The interval $[x,y-1]$ contains exactly $y-x$ integers.  
Consider the interval $I = [ g(x)  \log \theta , g(y) \log \theta)$. 
Because the length of $I$ is greater than $y-x+1$, it follows that $I$ contains at least $y-x+1$ integers $n$; that is, there are at least $y-x+1$ integers $n$ such that $x \leq M_{\theta}(n) \leq y-1$.  
By the pigeonhole principle, at least one of the intervals  
$\left[ g(x + i-1)\log \theta, g(x + i) \log \theta \right)$ 
with $i = 1,\ldots, y-x$ contains two integers, and so there exist integers $n$ and $n+1$ such that 
$M_{\theta}(n) = M_{\theta}(n+1)=i $.  
It follows that  if $\theta > e$, then the function $M_{\theta}(n)$ is not eventually strictly increasing.
This completes the proof.
$\blacksquare$ 

For example, if $\theta = 2$, then 
Theorem~\ref{FPR:theorem:StrictIncrease} and   
Corollary~\ref{FPR:corollary:BoundedGapsIneq}
imply that, for $n \geq 2$,  the function $M_2(n)$ is strictly increasing with bounded gaps, and that 
\[
1 \leq M_2(n+1) - M_2(n) < \frac{1}{\log 2} + \frac{3}{2} < 3.
\]
Thus,  $M_2(n+1) - M_2(n) \leq 2$ and so $\left(M_2(n+1) - M_2(n)-1\right)_{n=2}^{\infty}$ is a binary sequence, that is, a sequence of 0s and 1s.

\section{Explicit values and linear periodicity}  
\label{FPR:section:linearity}  

In Corollary~\ref{FPR:corollary:e} 
we proved that $M_e(n) = n-1$ for all $n \geq 2$.
This allows us to compute other explicit values of the function $M_{\theta}(n)$.  
For example, if $\ell \geq 2$ and $n \geq 1$, then 
\[
M_{e^{1/\ell}}(n) 
= \left[ \frac{1}{e^{1/\ell n}-1} \right] = M_e(\ell n) = \ell n-1.
\]

Let $k$ and $\ell$ be relatively prime positive integers, 
and let $\theta = e^{k/\ell}$.
Let $q\in \N$ satisfy $q > 1/(\ell \log 2)$.
If $n = kq$, then $n > \log \theta/\log 2$ and
\[
M_{e^{k/\ell}}(n) 
= \left[ \frac{1}{ e^{k/\ell n} -1} \right] 
= \left[ \frac{1}{ e^{1/\ell q} -1} \right] = M_{e}(\ell q) = \ell q -1.
\]
Because $n+k = k(q+1)$, we have 
\[
M_{e^{k/\ell}}(n+k) = M_{e^{k/\ell}}(n) + \ell.
\]
If $n = kq+r$, where $1 \leq r \leq k-1$, then
\[
\theta^{1/n} = e^{k/\ell n} = e^{(kq/n)/\ell  q}
= e^{(1 - r/n) / \ell q}
\]
and so
\[
M_{e^{k/\ell}}(n) = \left[  \frac{1}{e^{k/\ell n} - 1}  \right] 
= \left[  \frac{1}{  e^{(1 - r/n ) / \ell q} - 1}  \right] 
=  M_{   e^{(1 -r/n)  }}(\ell q).
\]

These results suggest examining the function $M_{\theta}(n)$ for numbers $\theta$ such that $\log \theta$ is rational.  
We compute $M_{\theta}(n)$ for $\log\theta = 2/3$,  $2/5$, $4/5$, $3/7$, and $5/7$,  and $1 \leq n \leq 90$.
We put a box around $M_{\theta}(N_0(\theta))$. 

\[
\boxed{\theta = e^{2/3}}
\]
\[
\begin{array}{ccccccccccccccc}
  \boxed{1} & 2   & 4   & 5   & 7   & 8   & 10 & 11 & 13 & 14 & 16 & 17 & 19 & 20 & 22 \\
23 & 25 & 26 & 28 & 29 & 31 & 32 & 34 & 35 & 37 & 38 & 40 & 41 & 43 & 44 \\ 
46 & 47 & 49 & 50 & 52 & 53 & 55 & 56 & 58 & 59 & 61 & 62 & 64 & 65 & 67 \\
68 & 70 & 71 & 73 & 74 & 76 &   77 &   79 &   80 &   82 &   83 &   85 &   86 & 88   & 89 \\
  91 &   92 &   94 &   95 &   97 &   98 & 100 & 101 & 103 & 104 & 106 & 107 & 109 & 110 & 112 \\
113 & 115 & 116 & 118 & 119 & 121 & 122 & 124 & 125 & 127 & 128 & 130 & 131 & 133 & 134 
 \end{array}
\]
\\
\[
\boxed{\theta = e^{2/5}}
\]
\[
\begin{array}{ccccccccccccccc}
\boxed{2}& 4& 7& 9& 12& 14& 17& 19& 22& 24& 27& 29& 32& 34& 37 \\ 
39& 42& 44& 47& 49& 52& 54& 57& 59& 62& 64& 67& 69& 72& 74 \\
77& 79& 82& 84& 87& 89& 92& 94& 97& 99& 102& 104& 107& 109& 112 \\
114& 117& 119& 122& 124& 127& 129& 132& 134& 137& 139& 142& 144& 147& 149 \\
152& 154& 157& 159& 162& 164& 167& 169& 172& 174& 177& 179& 182& 184& 187 \\ 
189& 192& 194& 197& 199& 202& 204& 207& 209& 212& 214& 217& 219& 222& 224
 \end{array}
\]
\\
\[
\boxed{\theta = e^{4/5}}
\]
\[
\begin{array}{ccccccccccccccc}
  4 &   \boxed{2} &   3 &  4  &  5  &  7  &   8 &   9 & 10 & 12 & 13 & 14 & 15 & 17 & 18 \\
19 & 20 & 22 & 23 & 24 & 25 & 27 & 28 & 29 & 30 & 32 & 33 & 34 & 35 & 37 \\
38 & 39 & 40 & 42 & 43 & 44 & 45 & 47 & 48 & 49 & 50 & 52 & 53 & 54 & 55 \\
57 & 58 & 59 & 60 & 62 & 63 & 64 & 65 & 67 & 68 & 69 & 70 & 72 & 73 & 74 \\
75 & 77 & 78 & 79 & 80 & 82 & 83 & 84 & 85 & 87 & 88 & 89 & 90 & 92 & 93 \\
94 & 95 & 97 & 98 & 99 & 100 & 102 & 103 & 104 & 105 & 107 & 108 & 109 & 110 & 112
 \end{array}
\]
\\
\[
\boxed{\theta = e^{3/7}}
\]
\[
\begin{array}{ccccccccccccccc}
 \boxed{1} &  4 &  6 &  8 &  11 &  13 &  15 &  18 &  20 &  22 &  25 &  27 &  29 &  32 &  34 \\
36 &  39 &  41 &  43 &  46 &  48 &  50 &  53 &  55 &  57 &  60 &  62 &  64 &  67 &  69 \\
71 &  74 &  76 &  78 &  81 &  83 &  85 &  88 &  90 &  92 &  95 &  97 &  99 &  102 &  104\\
106 &  109 &  111 &  113 &  116 &  118 &  120 &  123 &  125 &  127 &  130 &  132 &  134 &  137 & 139 \\
141& 144& 146& 148& 151& 153& 155& 158& 160& 162& 165& 167& 169& 172& 174 \\
176& 179& 181& 183& 186& 188& 190& 193& 195& 197& 200& 202& 204& 207& 209
\end{array}
\]
\\
\[
\boxed{\theta = e^{5/7}}
\]
\[
\begin{array}{ccccccccccccccc}
23 & \boxed{2} & 3 & 5 & 6 & 7 & 9 & 10 & 12 & 13 & 14 & 16 & 17 & 19 & 20\\ 
21 & 23 & 24 & 26 & 27 & 28 & 30 & 31 & 33 & 34 & 35 & 37 & 38 & 40 & 41\\ 
42 & 44 & 45 & 47 & 48 & 49 & 51 & 52 & 54 & 55 & 56 & 58 & 59 & 61 & 62 \\
63 & 65 & 66 & 68 & 69 & 70 & 72 & 73 & 75 & 76 & 77 & 79 & 80 & 82 & 83\\
84 & 86 & 87 & 89 & 90 & 91 & 93 & 94 & 96 & 97 & 98 & 100 & 101 & 103 & 104\\
 105 & 107 & 108 & 110 & 111 & 112 & 114 & 115 & 117 & 118 & 119 & 121 & 122 & 124 & 125
 \end{array}
\]
\\
We shall call an arithmetic function $f$ 
\emph{eventually linearly periodic} if there are positive integers $k$, $\ell$, and $n_0$ such that
\[
f(n+k) = f(n)+\ell
\]
for all $n \geq n_0$.  
We define the \emph{difference function} $\Delta(f)$ of an arithmetic function $f$ as follows:  $\Delta(f)(n) = f(n+1) - f(n)$.  The difference function is \emph{eventually periodic} if there are positive integers $k$ and $n_1$ such that $\Delta(f)(n+k) = \Delta(f)(n)$ for all $n \geq n_1$.  

For example,  consider the function $f$ whose sequence of values is $1,2,4,5,7,8,\ldots$, that is, 
\[
f(n) = 
\begin{cases}
3q-2 & \text{if $n = 2q-1$} \\
3q-1 & \text{if $n = 2q$.}
\end{cases}
\]
Then 
\[
f(n+2) = f(n)+3
\]
for all $n\in \N$, and so $f$ is eventually linearly periodic.
We can also write
\[
f(n) = \frac{3}{2}n + \chi(n \pmod{2}  )
\]
where 
\[
\chi(n \pmod{2} ) =
\begin{cases}
-1 & \text{if $n\equiv 0 \pmod{2}$ } \\
-\frac{1}{2}  & \text{if $n\equiv 1 \pmod{2}$.}
\end{cases}
\]
The sequence of values of the difference function $\Delta(f)$ is $1,2,1,2,1,2,\ldots$, and so  $\Delta(f)(n+2) = \Delta(f)(n)$ for $n \geq 1$, that is, $\Delta(f)$ is eventually periodic.  
Note that $f(n) = M_{e^{2/3}}(n)$ for $1 \leq n \leq 90$.

\bl    \label{FPR:lemma:almost}
Let $f$ be an arithmetic function.  
Let $k$, $\ell$, and $n_0$ be positive integers.
The following are equivalent:
\benum
\item
$f$ is eventually linearly periodic, and $f(n+k) = f(n)+\ell$ 
for all $n \geq n_0$;  
\item
there is a function $\chi$ defined on $\Z/k\Z$ such that
\[
f(n) = \frac{\ell}{k}n + \chi(n\pmod{k})
\]
for all  $n \geq n_0$;
\item
the difference function $\Delta(f)$ defined by 
$\Delta(f)(n) = f(n+1) - f(n)$ is eventually periodic, and 
$\Delta(f)(n+k) = \Delta(f)(n)$ for all $n \geq n_0$.
\eenum
\el

\textit{Proof.}
If $f$ is eventually linearly periodic, then there are positive integers $k$, $\ell$, and $n_0$ such that  $n_0 \equiv 0 \pmod{k}$ and $f(n+k) = f(n)+\ell$ for all  $n \geq n_0$.  
It follows that $f(n+qk) = f(n)+ q\ell$ for all $q \geq 0$ 
and $n \geq n_0$.  
For $r=0,1,\ldots, k-1$, we define
 $a_r = f(n_0+r)$ and
 \[
\chi(r \pmod k) = a_r - \frac{\ell}{k}  (n_0 + r).
 \]
If $n \geq n_0$, then there exist unique integers $q \in \N_0$ 
and $r\in \{0,1,\ldots, k-1\}$ 
such that $n = n_0 + qk + r$.   It follows that 
$n \equiv r \pmod{k}$ and 
\begin{align*}
f(n) & = f(n_0+qk + r) = f(n_0+r)+q\ell \\
& =  a_r + \left( \frac{n-n_0-r}{k} \right) \ell  \\
& = \frac{\ell}{k} n + \chi( n \pmod k).
\end{align*}
Conversely, this implies that
\begin{align*}
f(n+k) & = \frac{\ell}{k} (n+k) + \chi( n+ k \pmod k) \\
& = \frac{\ell}{k} n + \chi( n \pmod k) + \ell \\
& = f(n) + \ell
\end{align*}
for $n \geq n_0$, and so (1) and (2) are equivalent.  

Similarly,  if $f$ is eventually linearly periodic and, 
for all $n \geq n_0$, we have
\[
f(n+k) = f(n) + \ell
\]
then
\[
f(n+1+k) = f(n+1) + \ell
\]
and so 
\[
\Delta(f)(n+k) = f(n+1+k) - f(n+k) = f(n+1) - f(n) = \Delta(f)(n).
\]
Hence, $\Delta(f)$ is eventually periodic.  

Conversely, suppose that  $\Delta(f)(n+k) = \Delta(f)(n)$ for all $n \geq n_0$.  
Let $\ell = f(n_0 + k) - f(n_0).$
If $n \geq n_0$ and $ f(n + k) - f(n) = \ell$, then 
\begin{align*}
f(n+1 + k) - f(n+1) 
& = \sum_{i=1}^{k} \left( f(n+1+i) -  f(n+i) \right)  \\  
& = \sum_{i=1}^{k-1} \Delta(f)(n+i) + \Delta(f)(n+k)\\
& = \sum_{i=1}^{k-1} \Delta(f)(n+i) + \Delta(f)(n)\\
& = f(n+k) - f(n) \\
& = \ell.
\end{align*}
It follows by induction that $f(n+k) = f(n)+\ell$ for all $n \geq n_0$.
This proves that (1) and (3) are equivalent. 
$\blacksquare$ 

\bc    \label{FPR:corollary:Limit-kl}
Let $f$ be an arithmetic function.  
If $f$ is eventually linearly periodic, and if $k$, $\ell$, and $n_0$ are positive integers such that $f(n+k) = f(n)+\ell$ 
for all $n \geq n_0$, then
\[
\lim_{n\rightarrow\infty} \frac{f(n)}{n} = \frac{\ell}{k}.
\]
\ec

\textit{Proof.}
The function $\chi$ is bounded, and so
\[
\lim_{n\rightarrow\infty} \frac{f(n)}{n} 
= \lim_{n\rightarrow\infty}\left( \frac{\ell}{k} + \frac{\chi(n)}{n} \right) 
= \frac{\ell}{k}.
\]
$\blacksquare$ 

Let $\theta > 1$.  
The computational data suggest that $M_{\theta}(n)$ is eventually linearly periodic if there exist  positive integers $k$ and $\ell$ such that 
$\theta = e^{k/\ell}$.  
The data for $\theta =  e^{2/3}$ and $\theta =  e^{2/5}$ lead to the following explicit formula for $M_{\theta}(n)$ for  numbers of the form $\theta =  e^{2/\ell}$.

\bt     \label{FPR:theorem:2}
Let $\ell$ be an odd integer, $\ell \geq 3$, and let  $\theta = e^{2/\ell}$. 
For every positive integer $n$, 
\[
M_{ e^{2/\ell} }(n) = \frac{\ell}{2} n + \chi(n \pmod{2})
\]
where
\[
\chi(n  \pmod{2}) = 
\begin{cases}
- 1& \text{if $n \equiv 0 \pmod{2}$ } \\
- \frac{1}{2} & \text{if $n \equiv 1 \pmod{2}$.}
\end{cases}
\]  
\et

\textit{Proof.}
Let $x = M_{ e^{2/\ell}  }(n)$.  
Applying inequality~\eqref{FPR:MainIneq} with $\log\theta = 2/\ell$,
we obtain 
\[
x + \frac{1}{2} - \frac{1}{x}  
\leq \frac{\ell n}{2} < x + \frac{3}{2}.  
\]
If $n$ is even, then $\ell n/2$ is an integer.  
If $x \geq 3$, then 
\[
x < x + \frac{1}{6} \leq x + \frac{1}{2} - \frac{1}{x} \leq \frac{\ell n}{2}  < x+ 1 +\frac{1}{2}
\]
and so $\ell n/2 = x+1$; that is,
\[
M_{ e^{2/\ell} }(n) =  \frac{\ell}{2} n - 1.
\]
If $n$ is even and $x < 3$, then $\ell = 3$, $\theta = e^{2/3}$, $n=2$, and $x=2$.  
In this case we also have $M_{e^{2/3} }(2) = 2 = (\ell/2)n -1$.  

If $n$ is odd, then $(\ell n-1)/2$ is an integer.  
If $x \geq 2$, then  
\[
x-1 < x  - \frac{1}{2} \leq x  - \frac{1}{x} \leq \frac{\ell n-1}{2} < x + 1
\]
and so $x = (\ell n - 1)/2$; that is, 
\[
M_{\theta}(n) =   \frac{\ell}{2} n - \frac{1}{2}.
\]
If $n$ is odd and $x=1$, then $\ell = 3$, $\theta = e^{2/3}$, and $n=1$.  
In this case we also have $M_{e^{2/3} }(1) = 1 = (\ell/2)n - 1/2$.  
This completes the proof.
$\blacksquare$ 

A fundamental result of this paper is the following necessary and sufficient condition for the eventual linear periodicity of $M_{\theta}(n)$.

\bt   \label{FPR:theorem:LinearPeriodicity}
Let $\theta > 1$.  The arithmetic function $M_{\theta}(n)$ is eventually linearly periodic if and only if 
there exist positive integers  $k$ and $\ell$ such that 
$\theta = e^{k/\ell}$.  
\et

\textit{Proof.}
Let $k, \ell \in \N$ and $\theta = e^{k/\ell}$.  
By Theorem~\ref{FPR:theorem:monotonicity}, 
$\lim_{n\rightarrow \infty} M_{\theta}(n) = \infty$, and so there exists an integer $n_0 > k/(\ell \log 2)$ 
such that $M_{\theta}(n) > 2k$ for all integers $n \geq n_0$.  
Let $n \geq n_0$ and $M_{\theta}(n) = x$.  
Then $x > 2k$.
Applying  inequality~\eqref{FPR:MainIneq}
 to $\theta = e^{k/\ell}$ and $\log \theta = k/\ell$, we obtain 
\beq     \label{FPR:IneqAgain}
\frac{k}{\ell}\left( x+\frac{1}{2} -\frac{1}{x} \right) 
\leq n  < \frac{k}{\ell}\left( x+\frac{3}{2}  \right).
\eeq
The inequality on the left of ~\eqref{FPR:IneqAgain}
implies that 
\[
-1 < -\frac{2k}{x} \leq 2\ell n - 2kx - k.
\]
Because $2\ell n - 2kx - k$ is an integer, it follows that 
\[
0 \leq 2\ell n - 2kx - k
\]
and so
\[
\frac{k}{\ell}\left( x+\frac{1}{2}  \right) \leq n.
\]
Adding $k$ to each side of this inequality, we obtain
\[
\frac{k}{\ell} g(x+\ell) < \frac{k}{\ell}\left( x + \ell+\frac{1}{2}  \right) 
\leq n + k.
\]
Similarly, the inequality  on the right of ~\eqref{FPR:IneqAgain}
is equivalent to 
\[
2\ell n - 2kx < 3k.
\]
Because $2\ell n - 2kx$ and $3k$ are integers, we have 
\[
2\ell n - 2kx \leq 3k -1.
\]
It follows that 
\[
n \leq \frac{k}{\ell}\left( x+\frac{3}{2} - \frac{1}{2k} \right)
\]
and, because $x > 2k$, 
\begin{align*}
n + k & \leq \frac{k}{\ell} \left( x+ \ell + \frac{3}{2} - \frac{1}{2k} \right) \\
& < \frac{k}{\ell} \left( x+ \ell + \frac{3}{2} - \frac{1}{x+\ell + 1} \right) \\
& \leq \frac{k}{\ell} g(x+\ell + 1).
\end{align*}
The inequality
\[
 \frac{k}{\ell} g(x+\ell) <  n + k < \frac{k}{\ell} g(x+\ell + 1) 
\]
implies that $M_{\theta}(n+k) = x + \ell = M_{\theta}(n)+\ell.$  
Thus, the function $M_{\theta}(n)$ is eventually linearly periodic.  

Conversely, if $\theta > 1$ and  $M_{\theta}(n)$ is eventually linearly periodic, then there exist positive integers 
$k$, $\ell$ and  $n_0$ such that 
$M_{\theta}(n+k) = M_{\theta}(n) + \ell$ for all $n \geq n_0$.
It follows that 
\[
M_{\theta}(n+ qk) = M_{\theta}(n) + q\ell
\]
for every integer $n \geq n_0$ and every positive integer $q$.  
Applying inequality~\eqref{FPR:MainIneq} to $M_{\theta}(n) = x$ 
and $M_{\theta}(n+qk) = x + q\ell$, we obtain
\[
 \log \theta   \left( x+\frac{1}{2} - \frac{1}{x} \right)
 \leq n <   \log \theta   \left( x+\frac{3}{2} \right)
\]
and 
\[
\log \theta   \left( x+ q \ell + \frac{1}{2} - \frac{1}{x +  q \ell} \right) 
 \leq n +  q k   < \log \theta   \left( x+ q \ell + \frac{3}{2} \right).
\]
Combining these inequalities gives 
\[
\log \theta  \left( x+\frac{1}{2} - \frac{1}{x} \right) +  q k  
<  \log \theta   \left( x+ q \ell + \frac{3}{2} \right)
\]
and 
\[
\log \theta   \left( x+ q \ell + \frac{1}{2} - \frac{1}{x +  q \ell} \right) 
< \log \theta  \left( x + \frac{3}{2} \right) +  q k 
\]
and so
\[
\frac{ q k }{  q \ell + 1 + \frac{1}{x}  } < \log \theta 
< \frac{  q k }{  q \ell - 1 - \frac{1}{x +  q \ell}   }.
\]
Equivalently,
\[
\frac{ k }{   \ell + \frac{1}{q}\left(1 + \frac{1}{x} \right) } 
< \log \theta 
< \frac{  k }{ \ell - \frac{1}{q}\left(1 + \frac{1}{x +  q \ell} \right)  }.
\]
This inequality holds  for all positive integers $q$, and so 
$\log \theta =   k/\ell.$
This completes the proof.  
$\blacksquare$ 

\section{An algorithm for $M_{e^{k/\ell}}(n)$}

A second fundamental result of this paper is an algorithm to compute $M_{e^{k/\ell}}(n)$.

\bt     \label{FPR:theorem:algorithm}
Let $k$, $\ell$, and $n$  be positive integers.  
For each $r \in \{0,1,\ldots, k-1\}$ there exist unique  integers $u_r$ and $v_r$ such that 
\beq             \label{FPR:ur}
k + 2\ell r = 2k u_r + v_r
\eeq
and
\beq             \label{FPR:vr}
0 \leq v_r \leq 2k-1.
\eeq
Define $\chi:\Z/k\Z \rightarrow \Q$ as follows:  
If  $n \equiv r \pmod{k}$ for $r \in \{0,1,\ldots, k-1\}$, then 
\beq             \label{FPR:chi}
\chi(n \pmod{k}) = u_r - 1 - \frac{\ell r}{k} 
= - \frac{1}{2} - \frac{v_r}{2k}.
\eeq
If
\beq             \label{FPR:BeattyLowerBound}
n >  \max\left( \frac{k}{\ell \log 2}, \left( e^{k/\ell} - 1\right)(2k+1)  \right)
 \eeq
then
\beq             \label{FPR:Beatty}
M_{e^{k/\ell}}(n) = \frac{\ell}{k} n + \chi(n \pmod{k}) 
= \left[ \frac{ \ell}{k} n - \frac{1}{2} \right].
\eeq
\et

\textit{Proof.}
We begin with the observation that if $v_r$ and $x$ are integers such that $v_r < 2k < x$, then 
\[
\frac{v_r}{2k} + \frac{1}{x} 
\leq \frac{2k-1}{2k} + \frac{1}{x} 
= 1 - \left(\frac{1}{2k} -  \frac{1}{x} \right) < 1.
\]

Let $n$ satisfy inequality~\eqref{FPR:BeattyLowerBound}.
Let $x = M_{e^{k/\ell}}(n)$.  
Theorem~\ref{FPR:theorem:monotonicity} implies that $x  > 2k$.  
If $r \in \{0,1,\ldots, k-1\}$ and $n \equiv r \pmod{k}$, then there exists $q \in \N_0$ such that $n = kq+r.$
Inequality~\eqref{FPR:MainIneq} gives
\[
\frac{k}{\ell}  \left(x + \frac{1}{2}- \frac{1}{x}\right) 
\leq kq+r 
<  \frac{k}{\ell}  \left(x + \frac{3}{2} \right).
\]
Rearranging, we obtain 
\[
x + \frac{1}{2}- \frac{1}{x} \leq \ell q + \frac{\ell r}{k}
< x + \frac{3}{2}
\]
and so 
\begin{align*}
x + 1 -  \frac{k+2\ell r}{2k} - \frac{1}{x} 
& \leq \ell q  < x + 2 - \frac{k+2\ell r}{2k}.
\end{align*}
It follows from~\eqref{FPR:ur} and~\eqref{FPR:vr} that 
\[
x - 1 < x - \frac{v_r}{2k} - \frac{1}{x} \leq \ell q +u_r - 1
< x +1- \frac{v_r}{2k} \leq x+1.
\]
Because $\ell q +u_r - 1$ is an integer, it follows that 
\begin{align*}
M_{e^{k/\ell}}(n) & = x  =  \ell q + u_r - 1\\
& = \ell\left( \frac{n-r}{k} \right) + u_r - 1\\
& = \frac{\ell}{k}n + u_r - 1 - \frac{\ell r}{k} \\
& = \frac{\ell}{k}n + \frac{k + 2\ell r-v_r}{2k} - 1 - \frac{\ell r}{k} \\
& = \frac{\ell}{k}n - \frac{1}{2} - \frac{v_r}{2k} \\
& = \frac{\ell}{k}n + \chi(n\pmod{k}).
\end{align*}
Because $0 \leq v_r/2k < 1$, we have
\[
\frac{\ell}{k}n - \frac{1}{2} - 1 < x \leq \frac{\ell}{k}n - \frac{1}{2}
\]
and so
\[
x = \left[ \frac{\ell}{k}n - \frac{1}{2} \right].
\]
This completes the proof.
$\blacksquare$ 

We shall apply Theorem~\ref{FPR:theorem:algorithm} to compute
$M_{e^{3/7}}(n)$.  
With $k = 3$ and $\ell = 7$, we have
\begin{center}
\begin{tabular}{|c|c|c|c|}
\hline
$r$ & $u_r$ & $v_r$ & $\chi(r \pmod{3})$ \\  \hline
0 & 0 & 3 & -1 \\
1 & 2 & 5 & -4/3 \\
2 & 5 & 1 & -2/3 \\ \hline
\end{tabular}
\end{center}
and so 
\begin{align*}
M_{\theta}(n) 
& = \frac{7}{3} n + \chi(n\pmod{3}) \\
& = \begin{cases}
\frac{7}{3} n - 1 & \text{if $n \equiv 0 \pmod{3} $}\\
\frac{7}{3} n - \frac{4}{3} &  \text{if $n \equiv 1 \pmod{3} $}\\
\frac{7}{3} n - \frac{2}{3}  &  \text{if $n \equiv 2 \pmod{3} $}
\end{cases} \\
& = \begin{cases}
7q - 1 & \text{if $n = 3q$}\\
7q + 1 & \text{if $n = 3q+1$}\\
7q + 4 & \text{if $n = 3q+2$.}
\end{cases}
\end{align*}

\section{Problems and remarks}  \label{FPR:section:problems}  
\benum

\item
For $\theta > 1$ with $\log \theta $ irrational, find patterns in the sequence 
$\left( M_{\theta}(n) \right)_{n=1}^{\infty}$.
Is it possible to ``predict'' the value of $M_{\theta}(n)$?
How ``pseudo-random'' is the deterministic sequence 
$M_{\theta}(n)$?  
Describe the set of all sequences of the form $\left( M_{\theta}(n) \right)_{n=n_0}^{\infty}$ for $\theta > 1$.

\item
Let $m \in \N$.  
The sequence of integers $\mca = (a_n)_{n=1}^{\infty}$ is \emph{uniformly distributed modulo $m$} if
\[
\lim_{N\rightarrow \infty} \frac{1}{N}\sum_{n=1}^N 
\card \left( \{ n \in \{1,2,\ldots, N\} : a_n \equiv r \pmod{m}  \} \right) 
= \frac{1}{m}
\]
for all $r \in \{0,1,\ldots, m-1\}$.  
The sequence of integers $\mca = (a_n)_{n=1}^{\infty}$ is \emph{uniformly distributed} if  $\mca$ is \emph{uniformly distributed modulo $m$} for all $m \in \N$.  
Is the integer sequence $([1/ \{ \theta^{1/n} \}] )_{n=1}^{\infty}$ uniformly distributed  for almost all $\theta > 1$?

\item
Consider the binary sequence 
$(  M_2(n+1) - M_2(n)  -1)_{n=2}^{\infty}$.
How ``almost periodic'' is this seqence?

\item
The Bernoulli numbers are the coefficients in the Taylor series 
\[
\frac{x}{e^x -1} = \sum_{r=0}^{\infty}  \frac{B_r}{r! }x^r 
= 1 - \frac{x}{2} + \sum_{r=1}^{\infty} \frac{B_{2r}}{(2r)! }x^{2r}.
\]
This series converges for $|x| < 2\pi$.
Equivalently,
\[
\frac{1}{e^{x/n} -1} 
= \frac{n}{x} - \frac{1}{2} + \sum_{r=1}^{\infty} \frac{B_{2r}}{(2r)! } \frac{x^{2r-1}}{n^{2r-1}}
\]
Writing $\theta = e^x$, we obtain 
\[
\frac{1}{\theta^{1/n} -1} 
= \frac{n}{ \log \theta } - \frac{1}{2} + \sum_{r=1}^{\infty} \frac{B_{2r}}{(2r)! } \frac{(\log \theta)^{2r-1}}{n^{2r-1}}.
\]

\item
The \emph{Beatty sequence} associated with the real numbers 
$\alpha$ and $\beta$ is the sequence 
$([n\alpha + \beta])_{n=1}^{\infty}$.  
By Theorem~\ref{FPR:theorem:algorithm},
if $\theta > 1$ and $\log \theta$ is rational, then the sequence
$\left(M_{\theta}(n) \right)_{n=1}^{\infty}$ eventually coincides 
with the Beatty sequence $([n/\log \theta - 1/2])_{n=1}^{\infty}$.  
By Corollary~\ref{FPR:corollary:Beatty-like},
$M_{\theta}(n) = [n/\log \theta - 1/2]$ or $[n/\log \theta + 1/2]$
for all $n > \log\theta/\log 2$.  
Kevin O'Bryant and the author~\cite{nathobry12} have proved that there exist real numbers  $\theta > 1$ such that 
$M_{\theta}(n) = [n/\log \theta + 1/2]$ for infinitely many positive integers $n$, but they also proved that, for every $\theta > 1$,  $M_{\theta}(n) = [n/\log \theta - 1/2]$ for almost all positive integers $n$; that is, if $\theta > 1$, then the set 
\[
\left\{ n \in \N:  M_{\theta}(n) 
= \left[ \frac{n}{\log \theta}  + \frac{1}{2} \right] \right\}
\]
has asymptotic density 0.

\item
We could have considered the function 
$M'_{\theta}(n) =  \left[ 1/  \| {\theta}^{1/n} \|  \right]$ 
instead of $M_{\theta}(n) =  \left[ 1/  \{ {\theta}^{1/n} \}  \right]$.
However, for $\theta > 1$ and  $n \geq \log \theta/\log 3/2$, 
we have $1 < \theta^{1/n} \leq 3/2$ and so 
$ \{ {\theta}^{1/n} \} = \| {\theta} \| = \theta^{1/n} -1$.  
Thus,  the functions  
$M'_{\theta}(n)$ and $M_{\theta}(n)$ eventually coincide.

Let $0 < \theta < 1$ and let $\psi = \theta^{-1} > 1$.
Then
\[
\frac{1}{1 -  \theta^{1/n} } = \frac{1}{1 - \psi^{-1/n} } 
= \frac{\psi^{1/n}}{\psi^{1/n} - 1}  = \frac{1}{\psi^{1/n} - 1} + 1.
\]
If  $n \geq -\log \theta/\log  3/2$, then 
$2/3 \leq  \theta^{1/n} < 1$ and $1 < \psi^{1/n} \leq 3/2$;  
hence $\left\| \theta^{1/n} \right\| = 1 - \theta^{1/n}$ and
$\left\| \psi^{1/n} \right\| =  \psi^{1/n} - 1$.
It follows that 
\[
M'_{ \theta}(n) = \left[ \frac{1}{1 -  \theta^{1/n} } \right] 
= \left[ \frac{1}{\psi^{1/n} - 1} \right]  + 1 
= M'_{\psi}(n) + 1
\]
and it suffices to consider  $M'_{\theta}(n)$ only for $\theta > 1$.
Thus, there is no essential difference between the functions
$M_{\theta}(n)$and $M'_{\theta}(n)$.

\item
If $\mca = (\theta_n)_{n=1}^{\infty}$ is any sequence of real numbers, then we can examine the arithmetic function 
\[
M_{\mca}(n) = \left[\frac{1}{\{\theta_n\}}\right].
\]
If $1 < \theta_n < 2$, then 
\[
M_{\mca}(n) = \left[\frac{1}{\theta_n -1}\right].
\]
Consider, for example, the sequence 
$\mca = \left( n^{1/n} \right)_{n=1}^{\infty}$.  
For every integer $x \geq 2$, what is the smallest integer $n$ such that $M_{\mca}(n) = x$?
\eenum

\appendix

\section{Estimates for the exponential and logarithmic functions}  \label{FPR:appendix:exp}

This section contains the proofs of the estimates 
for the exponential and logarithmic functions that were used in Sections~\ref{FPR:section:growth}  
and~\ref{FPR:section:linearity}.

\bl    \label{FPR:lemma:exp(x)}
For all real numbers $x > 0$ and integers $n \geq 1$,
\beq  \label{FPR:ExpIneqLess}
\left(1+\frac{x}{n}\right)^n < \left(1+\frac{x}{n+1}\right)^{n+1} < e^x
\eeq
For all real numbers $x$ such that $0 < x \leq 1$ and for all integers $n \geq 2$,
\beq  \label{FPR:ExpIneqMore}
e^x < \left(1+\frac{x}{n}\right)^{n+1} < \left(1+\frac{x}{n-1}\right)^{n}.
\eeq
\el

\textit{Proof.}
Recall that 
\[
e^x = \lim_{n\rightarrow \infty} \left(1+\frac{x}{n}\right)^n
= \lim_{n\rightarrow \infty} \left(1+\frac{x}{n-1}\right)^n 
\]
for all real and complex numbers $x$. 
For $1 \leq i <  k \leq n$, we have
\[
0 < \frac{i}{n+1} < \frac{i}{n} < 1
\]
and so
\[
\prod_{i=1}^{k-1} \left( 1 -  \frac{i}{n}\right) 
< \prod_{i=1}^{k-1} \left( 1 -  \frac{i}{n+1} \right).
\]
If $x > 0$, then  the binomial theorem gives 
\begin{align*}
 \left(1 + \frac{x}{n}\right)^n 
 & =  \sum_{k=0}^n \prod_{i=1}^{k-1} \left(1 - \frac{i}{n} \right)  \frac{x^k}{k!} 
  \leq  \sum_{k=0}^n \prod_{i=1}^{k-1} \left(1 - \frac{i}{n+1} \right)  \frac{x^k}{k!} \\
 & <  \sum_{k=0}^{n+1} \prod_{i=1}^{k-1} \left(1 - \frac{i}{n+1} \right)  \frac{x^k}{k!} 
 =  \left(1 + \frac{x}{n+1}\right)^{n+1}
\end{align*}
and so 
\[
e^x = \sup\left\{   \left(1 + \frac{x}{n}\right)^n : n=1,2,\ldots \right\}.
\]
This proves~\eqref{FPR:ExpIneqLess}.

If $0 < x \leq 1$ and  $n \geq 2$, then $x < n/(n-1)$.  
Equivalently,
\[
1+\frac{x}{n} <  1 + \frac{nx }{ (n-1)(n+x) }.
\]
Again applying the binomial theorem, we obtain  
\begin{align*}
1+\frac{x}{n} & <  1 + \frac{nx }{ (n-1)(n+x) } 
<  \left( 1 + \frac{x }{ (n-1)(n+x) } \right)^{n}  \\
& = \left( \frac{n(n-1+x) }{ (n-1)(n+x) } \right)^{n} 
= \left( \frac{\frac{n-1+x}{n-1}}{ \frac{n+x}{n} } \right)^{n}   \\
& = \left( \frac{1+\frac{x}{n-1}}{ 1+\frac{x}{n} } \right)^{n}.
\end{align*}
It follows that 
\[
\left( 1+\frac{x}{n}  \right)^{n+1} <  \left( 1+\frac{x}{n-1}\right)^{n}
\]
and 
\[
e^x = \inf\left\{   \left(1 + \frac{x}{n-1}\right)^n : n=1,2,\ldots \right\}.
\]
This proves~\eqref{FPR:ExpIneqMore}.
$\blacksquare$ 

\bl      \label{FPR:lemma:g(x)ineq}
For $x \geq 1$, the function 
\[
g(x) = \frac{1}{\log\left(1 + \frac{1}{x}\right)}
\]
is positive and strictly increasing, and
\beq   \label{FPR:g1}
x+ \frac{1}{2} - \frac{1}{x} \leq  g(x) < x + \frac{1}{2}.
\eeq
For $y \geq x \geq 1$, 
\beq   \label{FPR:g2}
y - x - \frac{1}{y} < g(y) - g(x) < y - x + \frac{1}{x}.
\eeq
\el

\textit{Proof.}
For $x > 0$, the function $g(x)$ is positive, and is strictly increasing  because 
\[
g'(x) = \frac{1}{x(x+1)\log^2(1+1/x)} > 0.
\]

Let $t > -1$, and consider the function
\[
h(t) = \frac{4}{2+t} + \log(1+t).
\]
Because 
\[
h'(t) = -\frac{4}{(2+t)^2} + \frac{1}{1+t} = \frac{t^2}{(2+t)^2(1+t)} > 0
\]
it follows that $h(t)$ is strictly increasing and so $h(t) > h(0) = 2$ for $t > 0.$
Let $x > 0$ and $t=1/x$.  We obtain
\[
 \log\left(1 + \frac{1}{x}\right) > 2 - \frac{4x}{2x+1} = \frac{2}{2x+1}
\]  
and so 
\[
g(x) = \frac{1}{  \log\left(1 + \frac{1}{x}\right) } < x+\frac{1}{2}.
\]
This gives the upper bound for $g(x)$.

Let $0 < t < \delta < 1$.   
Using the Taylor polynomial of degree 1 for the function $\log(1+t)
$, we obtain a real number $u$ satisfying $0 < u < t$ such that 
\[
\log(1+t) = t -  \frac{t^2}{2(1+u)^2} <  t - \frac{t^2}{2 (1+\delta)^2}.
\]
It follows that 
\begin{align*}
\frac{1}{\log(1+t)} 
& > \frac{1}{ t - \frac{t^2}{2 (1+\delta)^2 } } 
 = \frac{1}{t} \left( \frac{1}{1-\frac{t}{2(1+\delta)^2} }  \right)  \\
& > \frac{1}{t} \left( 1 + \frac{t}{2(1+\delta)^2} \right)  
 = \frac{1}{t}  + \frac{1}{2(1+\delta)^2} \\
 & > \frac{1}{t}  + \frac{ (1-\delta)^2}{2} 
 > \frac{1}{t}  + \frac{1}{2} - \delta.
\end{align*}
Because this inequality is true for all $\delta > t$, we have
\[
\frac{1}{\log(1+t)} \geq \frac{1}{t}  + \frac{1}{2} - t.
\]
Replacing $t$ with $1/x$ gives the lower bound for $g(x)$. 
Inequality~\eqref{FPR:g2} is an immediate consequence
of~\eqref{FPR:g1}.
$\blacksquare$ 

\paragraph{Acknowledgements.}
I wish to thank Dakota Blair, Dick Bumby, and Kevin O'Bryant for helpful comments and discussions about this paper.

\def\cprime{$'$} \def\cprime{$'$} \def\cprime{$'$}
\providecommand{\bysame}{\leavevmode\hbox to3em{\hrulefill}\thinspace}
\providecommand{\MR}{\relax\ifhmode\unskip\space\fi MR }
\providecommand{\MRhref}[2]{%
  \href{http://www.ams.org/mathscinet-getitem?mr=#1}{#2}
}
\providecommand{\href}[2]{#2}

\bigskip

\noindent\textbf{Mel Nathanson} received a B.A. in philosophy 
from the University of Pennsylvania, started graduate school 
in biophysics at Harvard University and finished with a Ph.D. 
in mathematics from the University of Rochester.  
Intersecting  various mathematical subcultures, 
he studied with I. M. Gel'fand at Moscow State University, 
was Assistant to Andr\'e Weil at the Institute for Advanced Study, 
and wrote numerous papers with Paul Erd\H os.  
He has been on the faculty of Southern Illinois University at Carbondale, 
Rutgers University, and CUNY, and has had visiting appointments 
at Harvard, Princeton, Tel Aviv, and the Institute for Advanced Study.

\noindent\textit{Department of Mathematics, Lehman College (CUNY), Bronx, NY 10468 and CUNY Graduate Center, New York, NY 10016 \\
melvyn.nathanson@lehman.cuny.edu}

\end{document}